\input pstricks\input xy \xyoption{all}

\def\P {{\Bbb P}}
\newwrite\num
%
\output={\if N\header\headline={\hfill}\fi
\plainoutput\global\let\header=Y}
\magnification\magstep1
\tolerance = 500
\hsize=14.4true cm
\vsize=22.5true cm
\parindent=6true mm\overfullrule=2pt
\newcount\kapnum \kapnum=0
\newcount\parnum \parnum=0
\newcount\procnum \procnum=0
\newcount\nicknum \nicknum=1
\font\ninett=cmtt9

\font\ninebf=cmbx9

\font\sixbf=cmbx6
\font\ninesl=cmsl9

\font\nineit=cmti9

\font\ninerm=cmr9

\font\sixrm=cmr6
\font\ninei=cmmi9
\font\eighti=cmmi8
\font\sixi=cmmi6
\skewchar\ninei='177 \skewchar\eighti='177 \skewchar\sixi='177
\font\ninesy=cmsy9
\font\eightsy=cmsy8
\font\sixsy=cmsy6
\skewchar\ninesy='60 \skewchar\eightsy='60 \skewchar\sixsy='60
\font\titelfont=cmr10 scaled 1440
\font\paragratit=cmbx10 scaled 1200

\font\name=cmcsc10
\font\emph=cmbxti10

\font\tenmsbm=msbm10
\font\sevenmsbm=msbm7
%

%
\font\got=eufm10
\font\Got=eufm7
\font\teneufm=eufm10
\font\seveneufm=eufm7
\font\fiveeufm=eufm5
\newfam\eufmfam
\textfont\eufmfam=\teneufm
\scriptfont\eufmfam=\seveneufm
\scriptscriptfont\eufmfam=\fiveeufm

\font\tenmsam=msam10
\font\sevenmsam=msam7
\font\fivemsam=msam5
\newfam\msamfam
\textfont\msamfam=\tenmsam
\scriptfont\msamfam=\sevenmsam
\scriptscriptfont\msamfam=\fivemsam
\font\tenmsbm=msbm10
\font\sevenmsbm=msbm7
\font\fivemsbm=msbm5
\newfam\msbmfam
\textfont\msbmfam=\tenmsbm
\scriptfont\msbmfam=\sevenmsbm
\scriptscriptfont\msbmfam=\fivemsbm
\def\Bbb#1{{\fam\msbmfam\relax#1}}
\def\cz{{\kern0.4pt\Bbb C\kern0.7pt}
}
\def\ez{{\kern0.4pt\Bbb E\kern0.7pt}
}
\def\fz{{\kern0.4pt\Bbb F\kern0.3pt}}
\def\gz{{\kern0.4pt\Bbb Z\kern0.7pt}}
\def\hz{{\kern0.4pt\Bbb H\kern0.7pt}
}
\def\kz{{\kern0.4pt\Bbb K\kern0.7pt}
}
\def\nz{{\kern0.4pt\Bbb N\kern0.7pt}
}
\def\oz{{\kern0.4pt\Bbb O\kern0.7pt}
}
\def\rz{{\kern0.4pt\Bbb R\kern0.7pt}
}
\def\sz{{\kern0.4pt\Bbb S\kern0.7pt}
}
\def\pz{{\kern0.4pt\Bbb P\kern0.7pt}
}
\def\qz{{\kern0.4pt\Bbb Q\kern0.7pt}
}
\newskip\ttglue
\def\ninepoint{\def\rm{\fam0\ninerm}%
  \textfont0=\ninerm \scriptfont0=\sixrm \scriptscriptfont0=\fiverm
  \textfont1=\ninei \scriptfont1=\sixi \scriptscriptfont1=\fivei
  \textfont2=\ninesy \scriptfont2=\sixsy \scriptscriptfont2=\fivesy
  \textfont3=\tenex \scriptfont3=\tenex \scriptscriptfont3=\tenex
  \def\it{\fam\itfam\nineit}%
  \textfont\itfam=\nineit
  \def\sl{\fam\slfam\ninesl}%
  \textfont\slfam=\ninesl
  \def\bf{\fam\bffam\ninebf}%
  \textfont\bffam=\ninebf \scriptfont\bffam=\sixbf
   \scriptscriptfont\bffam=\fivebf
  \def\tt{\fam\ttfam\ninett}%
  \textfont\ttfam=\ninett
  \tt \ttglue=.5em plus.25em minus.15em
  \normalbaselineskip=11pt
  \font\name=cmcsc9
  \let\sc=\sevenrm
  \let\big=\ninebig
  \setbox\strutbox=\hbox{\vrule height8pt depth3pt width0pt}%
  \normalbaselines\rm
  \def\sl{\it}}

\headline={\ifodd\pageno\rightheadline\else\leftheadline\fi}
\def\rightheadline{\ninepoint Paragraphen"uberschrift\hfill\folio}
\def\leftheadline{\ninepoint\folio\hfill Chapter"uberschrift}
\let\header=Y
\def\titel#1{\need 9cm \vskip 2truecm
\parnum=0\global\advance \kapnum by 1
{\baselineskip=16pt\lineskip=16pt\rightskip0pt
plus4em\spaceskip.3333em\xspaceskip.5em\pretolerance=10000\noindent
\titelfont Chapter \uppercase\expandafter{\romannumeral\kapnum}.
#1\vskip2true cm}\def\leftheadline{\ninepoint
\folio\hfill Chapter \uppercase\expandafter{\romannumeral\kapnum}.
#1}\let\header=N
}
\def\Titel#1{\need 9cm \vskip 2truecm
\global\advance \kapnum by 1
{\baselineskip=16pt\lineskip=16pt\rightskip0pt
plus4em\spaceskip.3333em\xspaceskip.5em\pretolerance=10000\noindent
\titelfont\uppercase\expandafter{\romannumeral\kapnum}.
#1\vskip2true cm}\def\leftheadline{\ninepoint
\folio\hfill\uppercase\expandafter{\romannumeral\kapnum}.
#1}\let\header=N
}
\def\need#1cm {\par\dimen0=\pagetotal\ifdim\dimen0<\vsize
\global\advance\dimen0by#1 true cm
\ifdim\dimen0>\vsize\vfil\eject\noindent\fi\fi}
\def\neupara#1{\par\penalty-2000
\procnum=0\global\advance\parnum by 1
\vskip1cm\noindent{\paragratit \the\parnum. #1}%
\def\rightheadline{\ninepoint\S\the\parnum.\ #1\hfill \folio}%
\vskip 8mm\noindent}
\def\Proclaim #1 #2\finishproclaim {\bigbreak\noindent
{\bf#1\unskip{}. }{\it#2}\medbreak\noindent}
%
\gdef\proclaim #1 #2 #3\finishproclaim {\bigbreak\noindent%
\global\advance\procnum by 1
{%
{\relax\ifodd \nicknum
\hbox to 0pt{\vrule depth 0pt height0pt width\hsize
   \quad \ninett#3\hss}\else {}\fi}%
\bf\the\parnum.\the\procnum\ #1\unskip{}. }
{\it#2}
\immediate\write\num{\string\def
 \expandafter\string\csname#3\endcsname
 {\the\parnum.\the\procnum}}
\medbreak\noindent}
\newcount\stunde \newcount\minute \newcount\hilfsvar
\def\uhrzeit{
    \stunde=\the\time \divide \stunde by 60
    \minute=\the\time
    \hilfsvar=\stunde \multiply \hilfsvar by 60
    \advance \minute by -\hilfsvar
    \ifnum\the\stunde<10
    \ifnum\the\minute<10
    0\the\stunde:0\the\minute~Uhr
    \else
    0\the\stunde:\the\minute~Uhr
    \fi
    \else
    \ifnum\the\minute<10
    \the\stunde:0\the\minute~Uhr
    \else
    \the\stunde:\the\minute~Uhr
    \fi
    \fi
    }

 \def\calH{{\cal H}}
 \def\calJ{{\cal J}}
 \def\calL{{\cal L}}
 
\def\calO{{\cal O}} 
 
 \def\calT{{\cal T}}

\def\gota{\hbox{\got a}}

\def\gotm{\hbox{\got m}} 
\def\gotn{\hbox{\got n}}

\def\Gotm{\hbox{\Got m}}

\def\dim{\mathop{\rm dim}\nolimits}

\def\GL{\mathop{\rm GL}\nolimits}

\def\im{\mathop{\rm Im}\nolimits}

\def\kernel{\mathop{\rm kernel}\nolimits}

\def\mod{\mathop{\rm mod}\nolimits}

\def\proj{\mathop{\rm proj}\nolimits}
\def\rank{\mathop{\rm rank}\nolimits}

\def\Sp{\mathop{\rm Sp}\nolimits}

\def\boxit#1{
  \vbox{\hrule\hbox{\vrule\kern6pt
  \vbox{\kern8pt#1\kern8pt}\kern6pt\vrule}\hrule}}
\def\Boxit#1{
  \vbox{\hrule\hbox{\vrule\kern2pt
  \vbox{\kern2pt#1\kern2pt}\kern2pt\vrule}\hrule}}

\def\smallni{\smallskip\noindent }
\def\medni{\medskip\noindent }

\def\lo{\longrightarrow}

\def\loma{\longmapsto}
\def\betr#1{\vert#1\vert}

\def\imag{{\rm i}}
\def\pii{\pi {\rm i}}

\def\square{\hbox{\hbox to 0pt{$\sqcup$\hss}\hbox{$\sqcap$}}}
\def\qed{\ifmmode\square\else{\unskip\nobreak\hfil
\penalty50\hskip3em\null\nobreak\hfil\square
\parfillskip=0pt\finalhyphendemerits=0\endgraf}\fi}
\def\pn{\the\parnum.\the\procnum}
\def\downmapsto{{\buildrel
        {\vbox{\hbox{\hskip.2pt$\scriptstyle-$}}}
        \over{\raise7pt\vbox{\vskip-4pt\hbox{$\textstyle\downarrow$}}}}}

\def\DefJa{1.1}

\def\OZm{1.3}

\def\NatHI{1.6}

\def\OrbGz{2.2}

\def\TdetL{2.4}
\def\Tcopr{2.5}
\def\NewL{2.6}
\def\ProOpt{3.1}
\def\ProOptz{4.1}

\def\lamDet{4.3}
\def\Etdd{4.4}
\def\EtddV{4.5}

\def\AdSu{5.3}
\def\Labc{5.4}
\def\MT{5.5}

\def\WpAll{6.5}

\nopagenumbers
\immediate\newwrite\num
\nicknum=0  
\let\header=N
\def\tr{{\hbox{\rm tr}}}

\immediate\openout\num=theta-normal.num
\immediate\newwrite\num\immediate\openout\num=theta-normal.num
\def\RAND#1{\vskip0pt\hbox to 0mm{\hss\vtop to 0pt{%
  \raggedright\ninepoint\parindent=0pt%
  \baselineskip=1pt\hsize=2cm #1\vss}}\noindent}
\noindent
\centerline{\titelfont On the variety associated to the ring of}%
\vskip2mm
\centerline{\titelfont theta constants in genus 3}%
\def\leftheadline{\ninepoint\folio\hfill
Theta constants in genus 3}%
\def\rightheadline{\ninepoint Introduction\hfill \folio}%
\headline={\ifodd\pageno\rightheadline\else\leftheadline\fi}
\vskip 1.5cm
\leftline{\it \hbox to 6cm{Eberhard Freitag\hss}
Riccardo Salvati
Manni  }
  \leftline {\it  \hbox to 6cm{Mathematisches Institut\hss}
Dipartimento di Matematica, }
\leftline {\it  \hbox to 6cm{Im Neuenheimer Feld 288\hss}
Piazzale Aldo Moro, 2}
\leftline {\it  \hbox to 6cm{D69120 Heidelberg\hss}
 I-00185 Roma, Italy. }
\leftline {\tt \hbox to 6cm{freitag@mathi.uni-heidelberg.de\hss}
salvati@mat.uniroma1.it}
\vskip1cm\noindent%
\let\header=N%
\def\imag{{\rm i}}%
\def\transpose#1{\kern1pt{^t\kern-3pt#1}}%
\centerline{\vbox{\noindent\hsize=10cm
{\ninepoint{\bf Abstract}
\smallni
Due to fundamental results of Igusa [Ig1] and Mumford [Mu] the $N=2^{g-1}(2^g+1)$ even theta constants
define for each genus $g$ an injective holomorphic map of the Satake compactification
$X_g(4,8)=\overline{\calH_g/\Gamma_g[4,8]}$ into the projective space $P^{N-1}$. Moreover, this map is
biholomorphic onto the image outside the Satake boundary.
It is not biholomorphic on the whole in the cases $g\ge 6$ [Ig3].
Igusa also proved that in the cases $g\le 2$
this map biholomorphic onto the image [Ig2].
In this paper we extend this result to the case $g=3$. So we show that the theta map
$$X_3(4,8)\lo \P^{35}$$
is biholomorphic onto the image. This is equivalent to the statement that the
image is a normal subvariety of $\P^{35}$.}}}
\vskip5mm\noindent
{\paragratit Introduction}%
\medskip\noindent
The algebra $R(g,q)$  is generated by the theta constants
$$f_{a,q}=\sum_{n\in\gz^g}\exp\pii q Z[n+a/q].$$
Here $Z$ varies in the Siegel upper half space of genus $g$ and $a$ is a vector
in $\gz^g$. The series depends only on $\pm a$ mod $q$. We always assume that
$q$ is an even natural number. The functions $f_{a,q}$ are modular forms with respect
to the Igusa group $\Gamma_q[q,2q]$. In particular, the series $f_{a,q}/f_{b,q}$ are invariant under
$\Gamma_g[q,2q]$. The space of modular forms $[\Gamma_g[q,2q],r/2]$, $r\in\gz$, consists
of all holomorphic functions $f$ on the Siegel half plane $\calH_g$ such that $f/f_{a,q}^r$ are invariant, where
in the case $g=1$ the usual regularity condition
at the cusps has to be added.
The algebra of modular forms is
$$A(\Gamma_g[q,2q])=\bigoplus_{r\in\gz}[\Gamma_g[q,2q],r/2].$$
By a result of Baily,
the projective variety of the graded algebra $A(\Gamma_g[q,2q])$ can be identified,
as a complex space, with the Satake compactification of $\calH_g/\Gamma_g[q,2q]$,
$$\proj(A(\Gamma_g[q,2q]))=X_g(q,2q):=\overline{\calH_g/\Gamma_g[q,2q]}.$$
Due to basic theorems of Igusa [Ig1] and Mumford [Mu],
we have an everywhere regular, birational map
$$\overline{\calH_g/\Gamma_g[q,2q]}\lo \proj(R(q,g)).$$
This implies that $A(\Gamma_g[q,2q])$ is the normalization
of $R(q,g)$.
In the case $q=4$ this map is bijective and biholomorphic outside
the boundary. The case $q=2$ is exceptional. Here one knows
that the ring $R(g,2)$ is normal if $g\le 3$ [Ru].
Moreover,  $\proj(R(g,2))$ is not a normal variety when $g\geq 4$, [SM].
The ring  $R(g,4)$ is normal  if and only if $g\leq 2$ [Ig2, Ig5]. 
{Moreover the ideal of the relations is generated  by the   so called Riemann's relations.}
We shall
obtain the following main result.
\Proclaim
{Theorem}
{The map
$$\overline{\calH_3/\Gamma_3[4,8]}\lo \proj(R(3,4))$$
is biholomorphic.}
\finishproclaim
We mention that Igusa uses a slightly different setting. One can show that the ring
$R(g,q^2)$ can be generated the ``theta constants of first kind''
$$\sum_{n\>\hbox{\sevenrm integral}}\exp\pii(Z[n+a/q]+2b'(n+a/q)),\quad a,b\ \hbox{integral}.$$
In the case $g=3$, $q=2$, these are 36 different (up to sign) theta constants.
\smallskip
\noindent  {In a forthcoming  paper we shall consider the  projective  variety related to Riemann's relations in genus $g=3$.}
\neupara{Local rings of  modular varieties and their completion}%
We denote by
$$\calH_g=\{Z\in\cz^{g\times g};\quad Z=Z',\  \im Z\ \hbox{positive definite}\}$$
the Siegel upper half plane and by $\Sp(g,\gz)$ the Siegel modular group acting
on $\calH_g$ through $Z\mapsto (AZ+B)(CZ+D)^{-1}$.
Recall that the principal congruence  subgroup is defined as
$$\Gamma_g[q]=\kernel\bigl(\Sp(g,\gz)\lo\Sp(g,\gz/q\gz)\bigr)$$
and Igusa's subgroup as
$$\Gamma_g[q,2q]:=\bigl\{\;M\in\Gamma_g[q],\quad(CD')_0
\equiv(AB')_0\equiv0\;\mod2q\;\bigr\}.$$
Here $S_0$ denotes the column built of the diagonal of a square matrix $S$.
We generalize results from [FK] and [Kn].
We consider the Siegel modular variety $\calH_g/\Gamma_g[q,2q]$ and the
Satake compactification
$$X_g(q,2q)=\overline{\calH_g/\Gamma_g[q,2q]}.$$
For a decomposition $g=g_1+g_2$ we consider
the map
$$\calH_{g_1}\lo X_g(q,2q),\quad \tau\loma \lim_{t\to\infty}
\pmatrix{\tau&0\cr0&\imag t}.$$
We call the image of $\tau$ the standard boundary point
related to $\tau$.
The full Siegel modular group $\Sp(g,\gz)$ acts on $X_g(q,2q)$. Every boundary point is equivalent
to a standard boundary point. Hence we can restrict to study the standard boundary points.
We recall the description of the
analytic local ring of $X_g(q,2q)$ at such a point [Ig4].
\proclaim
{Definition}
{Let $U\subset\calH_{g_1}$ be an open subset and let $T$ be a semipositive integral
symmetric $g_2\times g_2$-matrix. The space $\calJ_T(U)$  consists of all
holomorphic functions $f:U\times\cz^{g_2\times g_1}\to\cz$ with the transformation property
$$\eqalign{
f(\tau,z+qh)&=f(\tau,z),\cr
f(\tau,z+qh\tau)&=\exp\{-\pii\tr(qT[h]\tau+2h'Tz)\}f(\tau,z)\quad \hbox{for}\ h\in\gz^{g_2\times g_1}.
\cr}$$
For a point $\tau_0\in\calH_{g_1}$ we define
$$\calJ_T(\tau_0)=\lim_{\longrightarrow} \calJ_T(U),$$
where $U$ runs through all open neighborhoods of $\tau_0$.}
DefJa%
\finishproclaim
In the case $T=0$ we have an everywhere holomorphic abelian function of $z$ which must be constant.
So we see
$$\calJ_0(\tau_0)=\calO_{\tau_0},$$
where $\calO_{\tau_0}$ denotes the local ring of the complex manifold $\calH_{g_1}$ at $\tau_0$.
In the case $q\ge 4$ we can identify $\calO_{\tau_0}$ with the local ring of
$\calH_{g_1}/\Gamma_{g_1}[q,2q]$ at the image of $\tau_0$, and we can consider $\calO_{\tau_0}$
as subring of the local ring of $X_g(q,2q)$ at the cusp related to $\tau_0$.
The spaces $\calJ_T(\tau_0)$ are modules over $\calO_{\tau_0}$,
moreover multiplication gives a map
$$\calJ_{T_1}(\tau_0)\otimes_{\calO_{\tau_0}}\calJ_{T_1}(\tau_0)\lo \calJ_{T_1+T_2}(\tau_0).$$
If we evaluate elements of the space $\calJ_T(\tau_0)$ at the point $\tau_0$ we get usual
spaces of theta functions.
\proclaim
{Definition}
{The space $J_T(\tau_0)$  consists of all
holomorphic functions $f:\cz^{g_2\times g_1}\to\cz$ with the transformation property
$$\eqalign{
f(z+qh)&=f(z),\cr
f(z+qh\tau_0)&=\exp\{-\pii\tr\;(q(T[h]\tau_0+2h'Tz)\}f(z)\quad \hbox{for integral}\ h.
\cr}$$}
DefJb
\finishproclaim
We have the evaluation map
$$\calJ_T(\tau_0)\lo J_T(\tau_0).$$
\proclaim
{Lemma}
{The $\calO_{\tau_0}$ modules
$\calJ_T(\tau_0)$ are finitely generated and free.
}
OZm%
\finishproclaim
{\it Proof.\/} Since the elements of $\calJ_T(\tau_0)$ are periodic in $z$, they admit a Fourier expansion
$$f(\tau,z)=\sum_{k\;\hbox{\sevenrm integral}}c_k\exp2\pii\tr( k'z)/q.$$
The Fourier coefficients are in $\calO_{\tau_0}$.
The second equation in Definition \DefJa\ gives
$$c_{k+qTh}=\exp(\pii\tr\,(qT[h]\tau+2k'h\tau)) c_k.$$
In the case that $T$ is invertible, one can prescribe the Fourier coefficients
$c_k$ for a system of representatives
mod $qT\gz^{g_2\times g_1}$ and then reconstruct $f$ as a linear combination of theta functions.
This shows that $\calJ_T(\tau_0)$ is free of finite rank. The case of a singular $T$ can be reduced
to the previous case in a standard way (taking a quotient by the nullspace of $T$).
\qed
\smallskip
The same argument gives generators of the vector space $J_T(\tau_0)$
and hence the following result.
\proclaim
{Lemma}
{The evaluation map
$$\calJ_T(\tau_0)\lo J_T(\tau_0)$$
is surjective.}
EvSur%
\finishproclaim
Another way to express this is
$$J_T(\tau_0)=\calJ_T(\tau_0)\otimes_{\calO_{\tau_0}}\cz.$$
Now we assume $q\ge 4$. Then the groups $\Gamma_g[q,2q]$ contain no element of
finite order besides the unit matrix. In this case
the analytic local ring $S^{\rm an}(\tau_0)$ of $X_g(q,2q)$ at the image of $\tau_0$ can be described as
the set of series
$$\eqalign{&
\sum_T a_T\exp(\pii\,\tr (TW)/q),\quad a_T\in \calJ_T(\tau_0),\cr&a_{T[U]}(\tau,z)=
a_T(\tau,zU')\quad\hbox{for}\ U\in\GL(g_2,\gz)[q],\cr}$$
where $T$ runs through all symmetric integral semipositive $g_2\times g_2$-matrices
and such that a certain convergence condition is satisfied [Ig4].
In the paper [Ig4] it has been shown that the ``Poincar\'e series''
$$H_{T,f}(\tau,z,W)=\sum f(\tau,WU')\exp(\pii\,\tr(WT[U])/q),\qquad f\in\calJ_T(\tau_0),$$
have this convergence property. The sum is taken over distinct $T[U]$
for $U\in\GL(g_2,\gz)[q]$. From the Supplement of Theorem 1 in [Ig4],
also the following result  follows.
\proclaim
{Proposition}
{The maximal ideal of the ring $S^{\rm an}(\tau_0)$ is generated by the Poincar\'e series
$H_{T,f}$ for non-zero $T$ and by the maximal ideal of the local local ring $\calO_{\tau_0}$.}
mIgP%
\finishproclaim
We introduce a filtration $\gotm_n$ on $S^{\rm an}(\tau_0)$.
For a semipositive integral $T$ we denote by $\lambda(T)$ the biggest number $k$ such that
$T$ can be written as $T=T_1+\cdots +T_k$ with non-zero integral and semipositive $T_i$.
In the case $T=0$ this is understood as $\lambda(T)=0$.
The associated filtration is
$$\gotn_n=\{P\in S^{\rm an}(\tau_0);\quad a_T=0\quad\hbox{for}\ \lambda(T)<n\}.$$
Then we define $\gotm_n$ to be the ideal generated by
$$\gotm(\calO_{\tau_0})^\mu\, \gotn_\nu,\quad \mu+\nu\ge n,$$
where $\gotm(\calO_{\tau_0})$ denotes the maximal ideal of $\calO_{\tau_0}$.
The ideal $\gotm=\gotm_1$ is the maximal ideal of $S^{\rm an}(\tau_0)$ and we have
$$\gotm_1\supset\gotm_2\supset\dots\quad\hbox{and}\quad \gotm_\mu\gotm_\nu\subset\gotm_{\mu+\nu}.$$
The Poincar\'e series $H_{T,f}$ is contained in $\gotm_n$ if either $\lambda(T)\ge n$ or if
$\lambda(T)<n$ and
$$f\in\gotm(\calO_{\tau_0})^{n-\lambda(T)}\calJ_T(\tau_0).$$
So we have
$$S^{\rm an}(\tau_0)/\gotm_n\cong\bigoplus_{\lambda(T)<n}\calJ_T(\tau_0)/\gotm(\calO_{\tau_0})^{n-\lambda(T)}\calJ_T(\tau_0).$$
We want to get rid of convergence conditions and therefore introduce a formal variant.
First we define
$$\hat \calJ_T(\tau_0)=\calJ_T(\tau_0)\otimes_{\calO_{\tau_0}}\hat\calO_{\tau_0}$$
where $\hat\calO_{\tau_0}$ denotes the completion of $\calO_{\tau_0}$.
Then we introduce the formal ring
$\hat S(\tau_0)$ that consists of all formal series
$$\eqalign{&
\sum_T a_T\exp(\pii\,\tr (TW)/q),\quad a_T\in \hat J_T(\tau_0),\cr&a_{T[U]}(\tau,z)=
a_T(\tau,zU')\quad\hbox{for}\ U\in\GL(g_2,\gz)[q].\cr}$$
The matrices $T$ run through all integral semipositive $g_2\times g_2$-matrices.
\smallskip
The ring $\hat S(\tau_0)$ is just the completion of $S^{\rm an}(\tau_0)$ with respect to the
filtration $(\gotm_n)$. We denote by $\bar S(\tau_0)$ the usual completion (by the powers of the maximal
ideal $\gotm$). From $\gotm^n\subset \gotm_n$ we obtain a natural homomorphism
$$\bar S(\tau_0)\lo \hat S(\tau_0).$$
\proclaim
{Theorem}
{The natural homomorphism
$$\bar S(\tau_0)\lo \hat S(\tau_0)$$
is an isomorphism.}
NatHI%
\finishproclaim
The case of the zero-dimensional boundary components has been treated
(in the more general context of arbitrary tube domains) by Kn\"oller [Kn] who refers to
[FK] where  the special case of the Hilbert modular group has been treated.
\smallni
{\it Proof of Theorem \NatHI.}
First we proof that the homomorphism is surjective.
Since $S^{\rm an}(\tau_0)/\gotm^k$ is a finite dimensional vector space we find for each
$k$ an $r$ such that
$$\gotm^k\cap \gotm_r=\gotm^k\cap \gotm_{r+1}=\cdots.$$
Therefore we can construct inductively a sequence of natural numbers $r_1<r_2<\cdots$
such that
$$\gotm_{r_1}\subset\gotm^2+\gotm_{r_2},\quad\gotm_{r_2}\subset\gotm^3+\gotm_{r_3},\quad
\gotm_{r_3}\subset\gotm^4+\gotm_{r_4},\quad\dots.$$
An arbitrary element $f\in \hat S(\tau_0)$ can be written in the form
$$f=f_{1}+f_{2}+\cdots,\quad f_{i}\in\gotm_{r_i}.$$
We construct inductively elements $g_j\in\gotm^j$, $a_j\in\gotm_{r_j}$ such that
$$f_1=g_1+a_2,\quad f_2+a_2=g_3+a_3,\quad f_3+a_3=g_4+a_4,\dots.$$
Then
$$f_1+f_2+\cdots f_{k-1}=g_1+g_2+\cdots g_k+a_k.$$
The series $g_1+g_2+\cdots$ converges in $\bar S(\tau_0)$. Its image in $\hat S(\tau_0)$ is $f$. This shows the
surjectivity. We now know that $\hat S(\tau_0)$ is noetherian too. To show injectivity it is enough that
the dimension of $\hat S(\tau_0)$ is greater or equal $\dim \bar S(\tau_0)=g(g+1)/2$.
Here we use the well-known result of commutative algebra that for every ideal $\gota$ in a noetherian ring $R$
that contains a non-zero divisor we have $\dim R>\dim R/\gota$.
The dimension of a local noetherian ring can be computed as the highest
coefficient of the Hilbert Samuel polynomial. Hence we must show that
$${\dim \hat S(\tau_0)/\hat{\gotm}^k\over k^{g(g+1)/2-1}},\qquad(\hat{\gotm}\ \hbox{maximal ideal of}\ \hat S(\tau_0),$$
is unbounded.
We define the ideals $\hat{\gotm}_k$ in $\hat S(\tau_0)$ in the same way
as the ideals $\gotm_k$ in $S^{\rm an}(\tau_0)$.
This means that we set
$$\hat{\gotn}_n=\{P\in \hat S(\tau_0);\quad a_T=0\quad\hbox{for}\ \lambda(T)<n\}$$
and $\hat{\gotm}_n$ to be the ideal generated by
$$\gotm(\hat\calO_{\tau_0})^\mu\, \hat{\gotn}_\nu,\quad \mu+\nu\ge n,$$
where $\gotm(\hat\calO_{\tau_0})$ denotes the maximal ideal of $\hat\calO_{\tau_0}$.
It is sufficient to show that
$${\dim \hat S(\tau_0)/\hat{\gotm}_k\over k^{g(g+1)/2-1}}$$
remains unbounded.
The description above by means of Poincar\'e series shows
$$\dim \hat S(\tau_0)/\hat{\gotm}_k=\dim S(\tau_0)/{\gotm}_k.$$
During the following estimates, $T$ always runs through a system of semipositive integral matrices
mod $\GL(g_2,\gz)[q])$ and $C_1,C_2\dots$ will denote suitable constants. We have
$$\eqalign{
\dim \hat S(\tau_0)/\hat{\gotm}_k&=\dim S(\tau_0)/{\gotm}_k\cr
&=\sum_{\nu+\lambda(T)=k}
\dim\calJ_T(\tau_0)/\gotm(\calO_{\tau_0})^\nu\calJ_T(\tau_0)\cr
&\ge C_1\sum_{\nu+\lambda(T)=k}\dim J_T(\tau_0)\nu^{g_1(g_1+1)/2}
.\cr}$$
We only keep $T$ which are invertible. The dimension of $J_T(\tau_0)$ then is then $\det(T)^{g_1}$ up to
a constant factor [Ig4]. We obtain
$$\ge C_2\sum_{\nu+\lambda(T)=k}(\det T)^{g_1}\nu^{g_1(g_1+1)/2}.$$
A trivial estimate states $\tr(T)\ge \lambda(T)$. We claim that also $(\det T)^{1/g_2}$ is greater or equal
than $\lambda(T)$ up to a constant factor. Since this statement is invariant under unimodular transformation,
it is sufficient to prove this for Minkowski reduced matrices. It follows from the standard inequalities
for Minkowski reduced matrices.
Therefore we get
$$\ge C_3\sum_{\nu+\lambda(T)=k}\lambda(T)^{g_1g_2}\nu^{g_1(g_1+1)/2}.$$
Now we restrict the summation to the range
$k/2\le\lambda(T)\le 3k/4$. Then $\nu\ge k/4$. Hence we get
$$\ge C_4 k^{g_1g_2+g_1/(g_1+1)/2}\>\#\{T; T\;\mod\; \GL(g_2,\gz)[q],\ k/4\le\lambda(T)\le 3k/4\}.$$
The asymptotic behaviour of the number of all $T$ with an upper bound for $\lambda$ has been determined
by Kn\"oller [Kn], Satz 2.3.1. This gives
$$\ge C_5 k^{g_1g_2+g_1/(g_1+1)/2}\cdot k^{g_2(g_2+1)/2}=C_5 k^{g(g+1)/2}.$$
This finishes the proof of Theorem \NatHI.\qed
\neupara{Optimal decompositions}%
We use the notation
$$\calT_g=\{\hbox{integral semipositive $g\times g$-matrices}\}.$$
The group $\GL(g,\gz)$ acts on $\calT_g$ through
$T\mapsto T[U]=U'TU$ from the right. In our context, matrices $T\in\calT_g$ of rank one
are important.
\smallskip
We call a non-zero element of $\calT_g$ {\it irreducible\/} if it cannot be written as sum of two
non zero elements of $\calT_g$.
We recall that
for a semipositive $T\in\calT_g$ we denote by $\lambda(T)$ the biggest number $k$ such that
$T$ can be written as $T=T_1+\cdots T_k$ with non-zero $T_i\in\calT_g$.
Notice that the irreducible elements $T$ are characterized by $\lambda(T)=1$
and that $\lambda(T)$ is invariant under unimodular transformations.
\proclaim
{Definition}
{Let T be a semipositive definite integral matrix. A decomposition
into irreducible integral matrices
$$T=T_1+\cdots+T_k,\quad\lambda(T)=k,$$
is called {\emph q-optimal} if all $T_i$ are of rank $1$ and  if for
arbitrary $U_1,\dots,U_k$ in $\GL(g,\gz)[q]$
one of the following two conditions holds.
\vskip1mm\item{\rm a)}
$\lambda(T_1[U_1]+\cdots+T_k[U_k])>k$.
\vskip1mm\item{\rm b)}
$\displaystyle T_1[U_1]+\cdots+T_k[U_k]\sim T\;\mod\;\GL(g,\gz)[q]$.
}
ProDe%
\finishproclaim
We will make use of the following two simple facts.
\vskip1mm
\item{1)} If $T=T_1+\cdots+ T_k$ is optimal then $T[U]=T_1[U]+\cdots+T_k[U]$ is optimal
for all $U\in\GL(g,\gz)$.
\vskip1mm
\item{2)} If $T=T_1+\cdots+ T_k$ is optimal then
$$\pmatrix{T&0\cr0&0}=\pmatrix{T_1&0\cr0&0}+\cdots+\pmatrix{T_k&0\cr0&0}$$
is optimal too.
\smallni
An integral matrix is called {\it primitive} if its entries are coprime.
Primitive semipositive matrices of rank 1 can be written as dyadic products
$$T=aa',\quad a\ \hbox{primitive column},$$
where $a$ is unique up to the sign.
The group $\GL(g,\gz)$
acts transitively on the set of all primitive columns. Hence it acts transitively on the set
of all primitive integral matrices of rank 1.
\proclaim
{Lemma}
{Let $q=2$ or $q=4$.
Two primitive semipositive integral matrices $T,S$ of rank one are equivalent mod $\GL(g,\gz)[q]$
if and only if they are congruent mod $q$ (i.e.~$T\equiv S\;\mod\;q$).
}
OrbGz%
\finishproclaim
{\it Proof.} We write $T,S$ in the form $T=aa'$, $S=bb'$. There must be an index $i$ such that
$a_i$ is odd. From $a_i^2\equiv b_i^2$ and $q=2,4$ we conclude $a_i\equiv \pm b_i$ mod $4$.
Since we can replace $b$ by $-b$ we can assume $a_i\equiv b_i$ mod $q$. Then $a_ia_j\equiv b_ib_j$
implies $a\equiv b$ mod $q$. Hence there exists a matrix $U\in\GL(g,\gz)[q]$ such that
$b=Ua$. This shows $S=T[U]$.\qed
\smallskip\proclaim
{Lemma}
{Let
$$T=\pmatrix{t_0&t_1\cr t_1&t_2},\quad 0\le t_1\le t_0,t_2,$$
be an integral semi positive matrix. Then
$$\lambda(T)=t_0+t_2-t_1.$$}
lamT%
\finishproclaim
{\it Proof.} The equality
$$T=t_1\pmatrix{1&1\cr 1&1}+(t_0-t_1)\pmatrix{1&0\cr0&0}+(t_2-t_1)\pmatrix{0&0\cr0&1}$$
shows $\lambda(T)\ge t_0+t_2-t_1$. We have to show the reverse inequality.
Let $T=T_1+\cdots+T_k$, $k=\lambda(T)$, where the $T_i$ are integral,
positive semidefinite and different from 0.
Consider the matrix
$$S=\pmatrix{1&-1/2\cr-1/2&1}.$$
It is positive definite. Obviously
$$t_0+t_2-t_1=\tr(TS)=\tr(ST_1)+\cdots+\tr(ST_k)\ge k=\lambda(T).$$
This implies $\lambda(T)\le t_0+t_2-t_1$.\qed
\smallskip
\proclaim
{Lemma}
{Assume that $T$ is a positive definite integral $2\times 2$-matrix.
Then
$$\lambda(T)\ge {3\over 2}\sqrt{\det T}.$$}
TdetL%
\finishproclaim
{\it Proof.} Since $\lambda$ and $\det$ are unimodular invariant, we can assume
that $T$ is Minkowski reduced ($0\le 2t_1\le t_0\le t_2$). Then
$$\lambda(T)=t_0+t_2-t_1\ge {3\over 4}(t_0+t_2)\ge{3\over 2}\sqrt{t_0t_2}
\ge{3\over 2}\sqrt{t_0t_2-t_1^2}.\eqno\square$$
\vskip2mm
\proclaim
{Lemma}
{Let
$$T=\pmatrix{t_1&t'\cr t&T_2},\quad T_2\in\calT_{g-1},$$
be an integral primitive semipositive $g\times g$-matrix of rank $1$. Assume
$$\lambda\pmatrix{t_1+1&t'\cr t&T_2}=2.$$
Then $T_2$ is primitive or zero.}
Tcopr%
\finishproclaim
{\it Proof.} After a suitable unimodular
transformation with a matrix of the form ${1\,0\choose 0\,U}$ we can assume that
$$T_2=\pmatrix{d&0\cr0&0}.$$
We have to show $d\le 1$.
We have
$$T=\pmatrix{t_1+1&s&0\cr s&d&0 \cr 0&0&0},$$
 hence
$$\det\pmatrix{t_1+1& s\cr s&d}=d$$
and the claim follows from Lemma \TdetL\qed
\proclaim
{Lemma}
{Let $T$ be a semipositive $g\times g$-matrix of rank one with coprime entries and let $U\in\GL(g,\gz)[2]$
such that $\lambda(T+T[U])=2$. Then $T[U]=T$.}
NewL%
\finishproclaim
{\it Proof.\/} This statement is invariant under $T\mapsto T[V]$ where $V\in\GL(g,\gz)$.
Hence we can assume that $T$ is the matrix with $t_{11}=1$ and zeros elsewhere.
Let
$$H=T[U]=\pmatrix{h_1+1&h'\cr h&H_2}.$$
The entries of $H_2$ are even and hence not coprime. From Lemma
\Tcopr\ follows that they are zero. This implies $H=T[U]=T$.\qed
\neupara{Degree two}%
We prove the existence of optimal decompositions in the case $g=2$.
\proclaim
{Proposition}
{In the cases $g=2$, $q$ arbitrary (even),
every semipositive integral matrix $T$ admits an optimal decomposition.}
ProOpt%
\finishproclaim
This Proposition  is invariant under unimodular transformation.
Hence it is enough to prove Proposition \ProOpt\ for invertible Minkowski-reduced $T$
(i.e.\ $0\le 2t_{12}\le t_{11}\le t_{22}$).
\vskip2mm\noindent
{\it Proof of Proposition \ProOpt.}
We can assume that $T$ is invertible and Minkowski reduced.
Then we claim that
$$T=(t_0-t_1)E_1+(t_2-t_1)E_2+t_1E_3\quad (k=\lambda(T)=t_0+t_2-t_1)$$
where
$$E_1=\pmatrix{1&0\cr0&0},\quad E_2=\pmatrix{1&0\cr0&0},\quad E_3=\pmatrix{1&1\cr1&1}$$
and
$$r_1=t_0-t_1,\quad r_2=t_2-t_1,\quad r_3=t_1.$$
is optimal (in both cases $q=2$ and $q=4$).
We write the decomposition of $T$ in the form $T_1+\cdots+T_k$ where
$T_i$ belong to $\{E_1,E_2,E_3\}$. We have to consider
$$\tilde T=T_1[U_1]+\cdots+T_k[U_k],\quad U_i\in\GL(2,\gz)[q].$$
We can assume that $\lambda(\tilde T)=\lambda(T)$. Then we have to show
that $T$ and $\tilde T$ are equivalent under the group $\GL(2,\gz)[q]$.
From Lemma \NewL\ we can assume that
$$T_i= T_k\Longrightarrow U_i=U_k.$$
Hence we can write
$$\tilde T=r_1E_1[U_1]+r_2E_2[U_2]+r_3E_3[U_3],\quad U_i\in\GL(2,\gz)[q].$$
We can assume that $U_1=E$ is the unit matrix. Since $T$ is Minkowski reduced,
$r_1=t_0-t_1$ and $r_2=t_2-t_1$ both are positive. From Lemma \Tcopr\ we see that $(E_2[U_2])_{11}\le 1$.
But since this expression is even, it must be zero. Then necessarily $E_2[U_2]=E_2$.
So we can assume $U_1=U_2=E$. In the case $r_3=0$ we are finished.
Otherwise, we can apply Lemma \Tcopr\ again to see that the diagonal elements
of $E_3[U_3]$ are $\le 1$. They are odd, hence both are 1.
So we get
$$E_3[U_3]=\pmatrix{1&\pm 1\cr \pm1&1}.$$
In case of the plus sign we are done. The minus sign only can
occur of $q\le 2$.
Then we can transform with
${1\;\;\;0\choose0\,-1}$ which is a matrix in $\GL(2,\gz)[2]$. This completes the proof
of Proposition \ProOpt.\qed
\neupara{Degree three}%
We prove the existence of optimal decompositions in the case $g=3$.
\proclaim
{Proposition}
{In the cases $g=3$, $q=2,4$,
every semipositive integral matrix $T$ admits an optimal decomposition.}
ProOptz%
\finishproclaim
A $3\times 3$-matrix symmetric positive definite real matrix $T$ is reduced in the
sense of Minkowski if
$$\eqalign{&t_{11}\le t_{22}\le t_{33},\cr&0\le 2t_{12}\le t_{11},\quad 0\le 2t_{23}\le t_{22},\quad
2\betr{t_{13}}\le t_{11},\cr& 2(t_{12}+t_{23}+\betr{t_{13}})\le t_{11}+t_{22}.\cr}$$
\proclaim
{Lemma}
{Let $T$ be a positive definite reduced integral $3\times3$-matrix. Then
$$\lambda(T)=\cases{t_{11}+t_{22}+t_{33}-t_{12}-t_{23}+t_{13}&if $t_{13}\le0$,\cr
t_{11}+t_{22}+t_{33}-t_{12}-t_{23}-t_{13}+\min(t_{12},t_{13},t_{23})&if $t_{13}>0$.\cr}$$}
LaTf%
\finishproclaim
{\it Proof.}
We introduce a basic system of matrices
{\ninepoint
$$\eqalign{&E_1=\pmatrix{1&0&0\cr0&0&0\cr0&0&0},\ E_2=\pmatrix{0&0&0\cr0&1&0\cr0&0&0},
\ E_3\pmatrix{0&0&0\cr0&0&0\cr0&0&1},\
E_4=\pmatrix{1&1&0\cr1&1&0\cr0&0&0},\cr&
E_5=\pmatrix{0&0&0\cr0&1&1\cr0&1&1},\ E_6=\pmatrix{1&0&1\cr0&0&0\cr1&0&1},
\ E_7=\pmatrix{1&1&1\cr1&1&1\cr1&1&1}.\cr}$$}%
It is a system of representatives of integral semipositive matrices of rank one with respect
to the
action $T\mapsto T[U]$ of the group $\GL(2,\gz)[2]$. We also introduce the modified matrix
$$E_6^-=\pmatrix{1&0&-1\cr0&0&0\cr-1&0&1}.$$%
In the case $t_{13}\le0$ we use the decomposition
$$\eqalign{
T=&(t_{11}-t_{12}+t_{13})E_1+(t_{22}-t_{12}-t_{23})E_2+(t_{33}+t_{13}-t_{23})E_3\cr+&
t_{12}E_4+t_{23}E_5-t_{13}E_6^-.\cr}$$
 We notice that $E_7$ does not occur in this decomposition.
Since the coefficients are nonnegative, we get $\lambda(T)\ge t_{11}+t_{22}+t_{33}-t_{12}-t_{23}+t_{13}$.
For the reverse inequality we use the positive matrix
$$S={1\over 2}\pmatrix{2&-1&1\cr-1&2&-1\cr1&-1&2}.$$
For each nonzero semipositive integral matrix $H$ the trace $\tr(SH)$ is a positive integer.
This implies $\tr(SH)\ge \lambda(H)$. In our case we get
$\lambda(T)\le \tr(ST)=t_{11}+t_{22}+t_{33}-t_{12}-t_{23}+t_{13}$.
This completes the proof in the first case $t_{13}\le 0$.
\smallskip
In the second case, $t_{13}>0$, we use a similar decomposition.
Setting  $m = \min \{t_{12},t_{23},t_{13}\}$, we take
$$\eqalign{&T=\cr&
(t_{11}-t_{12}-t_{13}+m)E_1+(t_{22}-t_{12}-t_{23}+m)E_2+
(t_{33}-t_{23}-t_{13}+m)E_3\cr&+(t_{12}-m)E_{4}+ (t_{23}-m)E_{5}+ (t_{13}-m)E_{6}
\cr&+m E_{7}\cr}$$
which shows
$\lambda(T)\ge  t_{11}+t_{22}+t_{33}-t_{12}-t_{13}-t_{23}+m$.
We observe that at least one of the  coefficients of
$E_4, E_5, E_6$  is 0.
To prove the reverse inequality one uses $\tr(ST)\ge\lambda(T)$ for one of the
following three matrices
$$S=\pmatrix{2&0&-1\cr0&2&-1\cr -1&-1&2}\ \hbox{or}\ \pmatrix{2&-1&-1\cr-1&2&0\cr -1&0&2}
\ \hbox{or}\ \pmatrix{2&-1&0\cr-1&2&-1\cr 0&-1&2}$$
depending on whether $m$ is $t_{12}$ or $t_{23},$ or $t_{13}$.\qed
\proclaim
{Lemma}
{Let $T$ be an integral positive definite $3\times 3$-matrix. Then
$$\lambda(T)^3\ge 8\det T.$$}
lamDet%
\finishproclaim
{\it Proof.} We can assume that $T$ is reduced. From the inequality
$$2(t_{12}+t_{23}+\betr{t_{13})}\le t_{11}+t_{22}\le t_{11}+t_{33}\le t_{22}+t_{33}$$
together with the trivial inequality
$$\lambda(T)\ge (t_{11}+t_{22}+t_{33}-t_{12}-t_{23}-\betr{t_{13}})$$
we get
$$\lambda(T)\ge2{t_{11}+t_{22}+t_{33}\over 3}\ge 2(t_{11}t_{22}t_{33})^{1/3}.$$
The statement of the Lemma now follows from Hadamard's
inequality
$$t_{11}t_{22}t_{33}\ge\det T.\eqno\square$$
\proclaim
{Lemma}
{Let $T$ be a  matrix of rank~$\le 1$.
Then
$$\det\left(E_1+E_2+T\right)=t_{33}.$$
In addition, let $T$ be semipositive and integral. Then one of
the following two inequalities hold.
\vskip1mm
\item{\rm a)} $\lambda\left(E_1+E_2+T\right)>3$.
\item{\rm b)} $t_{33}\le 1$.}
Etdd%
\finishproclaim
{\it Proof.} The computation of the determinant
is easy. Hence we have to prove only the second statement.
We assume that $t_{33}>1$. Since it is a square, we obtain $t_{33}\ge4$.
But $t_{33}$ is the determinant of the matrix . Hence we get from Lemma \lamDet\ that
$\lambda^3\ge 32$. This gives $\lambda(T)>3$.\qed
\smallskip
We will apply several times not only Lemma \Etdd\
but also an obvious generalization.
Let $U\in\GL(3,\gz)$. Then one has for $\rank T\le 1$
$$\det\left(\pmatrix{E&0\cr0&0}[U]+T\right)=T[U^{-1}]_{33}.$$
\proclaim
{Corollary of Lemma \Etdd}
{Assume that $T$ is a semipositive integral matrix of rank one. Then
$$\eqalign{\lambda(E_1+E_2+T)= 3&\Longrightarrow  t_{33}\le 1,\cr
\lambda(E_1+E_3+T)= 3&\Longrightarrow  t_{22}\le 1,\cr
\lambda(E_1+E_4+T)= 3&\Longrightarrow  t_{33}\le 1,\cr
\lambda(E_1+E_5+T)= 3&\Longrightarrow  t_{22}-2t_{23}+t_{33}\le 1,\cr
\lambda(E_2+E_3+T)= 3&\Longrightarrow  t_{11}\le 1,\cr
\lambda(E_2+E_4+T)= 3&\Longrightarrow  t_{33}\le 1,\cr
\lambda(E_2+E_5+T)= 3&\Longrightarrow  t_{11}\le 1,\cr
\lambda(E_3+E_4+T)= 3&\Longrightarrow  t_{11}-2t_{12}+t_{22}\le 1,\cr
\lambda(E_3+E_5+T)= 3&\Longrightarrow  t_{11}\le 1,\cr
\lambda(E_4+E_5+T)= 3&\Longrightarrow  t_{11}+t_{22}+t_{33}+2t_{12}-2t_{23}+2t_{13}\le 1.\cr
}$$}
EtddV%
\finishproclaim
For the proof of Proposition \ProOptz\ we can assume that $T$ is positive definite.
The proposition is invariant under arbitrary unimodular transformation.
Hence we can assume that $T$ is reduced.
We have to differ between the two cases:
\smallni
{\bf Case A.} $t_{13}\le 0$.
\hfill\break
{\bf Case B.} $t_{13}>0$.
\smallni
We start with case A. We use the decomposition
$$\eqalign{
T=&r_1E_1+r_2E_2+r_3E_3+r_4E_4+r_5E_5+r_6E_6^-,\cr
&r_1=t_{11}+t_{13}-t_{12},\quad r_2=t_{22}-t_{12}-t_{23},\quad r_3=t_{33}+t_{13}-t_{23},\cr&
r_4=t_{12},\quad r_5=t_{23},\quad r_6=-t_{13}.\cr}$$
We will show that it is $q$-optimal in both cases $q=2$ and $q=4$.
The reduction inequalities read as
\medni
\halign{\qquad$#\le\>$\hfil&$#$\hfil&\quad$#\le\>$\hfil&$#$\hfil\cr
r_4&r_1+r_6,&r_1+r_6&r_2+r_5,\cr
r_6&r_1+r_4,&r_2+r_4&r_3+r_6,\cr
r_5&r_2+r_4,&r_5+r_6&r_1+r_2.\cr}
\medni
Since $T$ is positive definite we have also that the diagonal elements are positive,
in particular
$$r_1+r_4+r_6>0.$$
We also mention that at least two of the coefficients $r_1$, $r_2$, $r_3$ do not
vanish. More precisely we state.
\smallni
Only the following 4 cases are possible.
\vskip1mm
\item{1)} $r_1>0$, $r_2>0$, $r_3>0$.
\item{2)} $r_1>0$, $r_2>0$, $r_3=0$ and $r_6>0$, $r_4=0$, $r_1=r_2=r_5=r_6$.
\item{3)} $r_1>0$, $r_2=0$, $r_3>0$ and $r_6=0$, $r_1=r_4=r_5$.
\item{4)} $r_1=0$, $r_2>0$, $r_3>0$ and $r_6>0$, $r_4=r_6$.
\smallni
For the proof one has to discuss the three cases $r_i=0$ separately. We start with
\smallni
Case 1) There is nothing to prove.
\smallni
Case 2)  $r_3=0$: . Thus  $t_{33}-t_{23}+ t_{13}=0$.
Hence$t_{23}=-t_{13}= t_{33}/2$, thus  $r_6>0$
and  by the  basic inequalities
$$2(t_{12}+t_{33}\leq t_{11}+t_{22}\leq 2t_{33}.$$
Hence  $$(t_{12}=0;  t_{33}= t_{11}=t_{22} $$
This implies $r_4=0, r_1=r_2=r_5=r_6.$
\smallni
We observe that in this  case the  matrix  $T$
has the form
$$ T= \pmatrix{2a&0&-a \cr 0&2a&a\cr -a&a&2a}$$
Case 3 )and Case 4) can be  proved in  similar way , we just observe that the  corresponding matrices $T$
have the forms
$$ T= \pmatrix{2a&a&0 \cr a&2a&a\cr 0&a&c}\quad  T= \pmatrix{2a&a&-a \cr a&b&h\cr -a&h&c} $$
\smallskip
Now we will prove that the described decomposition
$$T=r_1E_1+\cdots r_5E_5+r_6E_6^-$$
is $q$-optimal in each of the 4 cases. As in the case $g=2$ we can apply Lemma~\NewL\ to formulate
Proposition \ProOptz\ as follows. Consider matrices $U_1,\dots,U_6\in\GL(3,\gz)[q]$ and
$$\tilde T=r_1E_1[U_1]+\cdots r_5E_5[U_5]+r_6E_6^-[U_6].$$
Assume $\lambda(T)=\lambda(\tilde T)$. Then $T\sim\tilde T$ mod $\GL(3,\gz)[q]$.
\smallni
{\it Proof of Proposition \ProOptz\ in the case A1.}
\smallni
Without loss of generality we can assume that in the decomposition of $\tilde T$ we have $U_1=E$.
Then we have $\lambda(E_1+E_2[U_2])=2$ since this sum is a partial sum of $\tilde T$.
Now Lemma \Tcopr\ shows that
$$E_2[U_2]=H=\pmatrix{*&*\cr *&H_2}\quad\hbox{where}\ H_2\ \hbox{is primitive}.$$
(The other case in Lemma \Tcopr, $H_2=0$, cannot arise since the first diagonal element
of $H_2$ is odd.) We have the freedom to act on $H$ with a matrix of the form
${1\,0\choose 0\,V}$ where $V\in\GL(2,\gz)[q]$ since this does not change $E_1$.
Thanks to Lemma \OrbGz\ we can replace $H_2$ by the matrix
${1\,0\choose 0\,0}$. Since $H$ has rank one, $h_{11}$ must be zero. The semidefinitness now
implies $H=E_2$. Hence we can assume now $U_1=U_2=E$. Now we use $\lambda(E_1+E_2+E_3[U_3])=3$.
Lemma \Etdd\ shows $E_3[U_3]_{33}=1$. (Zero is not possible since this element is odd.)
We still can apply transformations with matrices of $\GL(3,\gz)[2]$ if they fix $E_1$ and $E_2$.
Hence we can multiply simultaneously the third row (column) by a multiple of $q$
and add it to another row (column).
This allows us to assume
$$E_3[U_3]=\pmatrix{*&*&0\cr*&*&0\cr0&0&1}.$$
Since the rank is one we get $E_3[U_3]=E_3$. Hence we can assume $U_1=U_2=U_3=E$.
Next we apply Lemma \EtddV\  to show that all diagonal elements
of the matrices $E_i[U_i]$, $i>3$, are 0 or 1. This shows that
$$E_4[U_4]=E_4[D_4],\quad E_5[U_5]=E_5[D_5],\quad E_6^-[U_6]=E_6^-[D_6]$$
where $D_i$ are diagonal matrices in $\GL(3,\gz)$.
In the case $q=4$ we are finished since then the congruence mod 4 shows $D_i=E$.
So we can assume $q=2$,
The diagonal matrices
fix $E_1,E_2,E_3$. Hence we can assume first $D_4=E$ and then $D_5=E$. There remain two possibilities
$E_6^-[D_i]=E_6$ or $E_6^-[D_i]=E_6^-$. The second is what we want, hence it remains to discuss
$E_6^-[D_i]=E_6$. In this case we claim that one of the $r_4,r_5$ is zero. Otherwise
$E_4+E_5+E_6=E_1+E_2+E_3+E_7$
would be a partial sum of $\tilde T$ which is not possible. So assume $r_4=0$. Then there is a diagonal
matrix $D$ with the property $E_6[D]=E_6^-$ which does not change anything in the first five summands.
This finishes the proof of A1.
\smallni
{\it Proof of Proposition \ProOptz\ in the case A2.}
\smallni
The decomposition of $T$ reads as
$$T=r_1(E_1[U_1]+E_2[U_2]+E_5[U_5]+E_6^-[U_6]).$$
As in the case A1 it is no loss of generality to assume $U_1=U_2=U_5=E$. Let
$H=E_6^-[U_6]$. From $\lambda(E_1+E_2+H)=3$ and Lemma \EtddV\ follows $h_{33}=1$
and similarly from $\lambda(E_2+E_5+H)=3$ follows $h_{11}=1$. Since
$h_{11}h_{33}=h_{13}^2$ we have $h_{13}=\pm1$. But $t_{13}\le0$, hence $h_{13}=-1$.
The matrix $H$ is semidefinit of rank 1. Hence it is of the form
$$H=\pmatrix{1&a&-1\cr a&a^2&-a\cr-1&-a&1}.$$
Now we use $\lambda(E_1+E_5+H)=3$. Lemma \Etdd\ shows $(a+1)^2\le 1$. Since $a$ is even, we get
$a=0$ or $a=-2$. In the case $a=0$ we are done. The case $a=-2$ occurs only if
$q=2$. Then we can
apply the transformation
$$\pmatrix{-1&0&0\cr0&1&0\cr0&-2&-1}\in\GL(3,\gz)[2].$$
It fixes $E_1,E_2,E_5$ and sends $H$ to $E_6^-$. This finishes the
proof of A2.
\smallni
{\it Proof of Proposition \ProOptz\ in the case A3.}
\smallni
We have
$$T=r_1(E_1+E_4+E_5)+r_3E_3\ \hbox{and}\ \tilde T=r_1(E_1[U_1]
+E_4[U_4]+E_5[U_5])+r_3E_3[U_3].$$
Again we can assume $U_1=E$. Considering the partial sum $T_1+T_3[U_3]$ we can
reduce to $U_3=E$. Then, considering $E_1+E_3+E_4[U_4]$, we get $E_4[U_4]_{22}\le 1$.
It must be 1 since it is odd.
Now, applying to $E_4[U_4]$ a unimodular substitution from $\GL(3,\gz)[q]$
that fixes $E_1$ and $E_3$, we can get $U_4=E$. So we can assume
$$\tilde T=r_1(E_1+E_4+E_5[U_5])+r_3E_3.$$
Now we apply Lemma \EtddV\  to
$$(E_1+E_3)+E_5[U_5],\quad (E_1+E_4)+E_5[U_5]$$
to obtain that $E_5[U_5]_{22}=1$ and $E_5[U_5]_{33}=0$. This means
$$E_5[U_5]=\pmatrix{a^2&0&a\cr 0&0&0\cr a&0&1}.$$
We have $a\equiv0$ mod $q$. We transform with the matrix from $\GL(2,\gz)[q]$.
$$\pmatrix{1&0&0\cr0&1&0\cr1-a&0&1}.$$
This transforms $\tilde T$  to $T$. This completes the proof of A3.
\smallni
{\it Proof of Proposition \ProOptz\ in the case A4.}
\smallni
We have
$$\eqalign{
T&=r_2E_2+r_3E_3+r_4E_4+r_5E_5+r_4E_6^-,\quad r_2,r_3,r_4>0,\cr
\tilde T &=r_2E_2[U_2]+r_3E_3[U_3]+r_4E_4[U_4]+r_5E_5[U_5]+r_4E_6^-[U_6].\cr}$$
Similar to the previous cases we can assume $U_2=U_3=E$.
Since $E_2+E_3+E_4[U_4]$ is optimal, we get $E_4[U_4]_{11}=1$. We can transform $\tilde T$ by
a matrix from $\GL(3,\gz)[q]$ that fixes $E_2,E_3$. This means that we can multiply the first
row (resp.\ column) of $E_4[U4]$ by a factor which is a multiple of $q$ and add it to the second
(or third row). In this way we can get
$$E_4[U_4]=\pmatrix{1&1&0\cr 1&*&*\cr 0&*&*}.$$
Since it is matrix of rank one, we then have $E_4[U_4]=E_4$. So we can assume $U_4=E$.
Now we assume $r_5>0$.
Then we can apply Lemma \Etdd\ to $E_2+E_3+E_5[U_5]$ to obtain
$E_5[U_5]_{11}=0$ (it is even) which implies
$$E_5[U_5]=\pmatrix{0&0\cr 0&H}.$$
The $2\times2$-matrix $H$ is primitive, semidefinit of rank one and its entries
are $\equiv 1$ mod $q$. Now Lemma \OrbGz\
implies that it is of the form
$$H=\pmatrix{1&1\cr1&1}[U],\quad U\in\GL(2,\gz)[q].$$
So we can assume
$r_5=0$ or $U_5=E$. It remains to treat $E_6^-[U_6]$. Since
$E_2+E_3+E_6^-[U_6]$ are optimal, we get $E_6[U_6]_{11}=1$.
Since $E_2+E_4+E_6^-[U_6]$ is optimal, we get $E_6[U_6]_{33}=1$.
Since this matrix is symmetric and of rank 1, it is of the form
$$E_6^-[U_6]=\pmatrix{1&a&\pm 1\cr a&a^2&\pm a\cr
\pm1&\pm a&1},\quad a\equiv 0\;\mod\; q.$$
In the case $q=0$ the minus sign must be there.
The case $r_5=0$ can be transformed to the case A2. (Interchange the first and the
third row and column). Hence we can assume $r_5>0$. Then we can consider
$E_4+E_5+E_6^-[U_6]$ which is optimal. Lemma \EtddV\ gives
$a^2\le 1$ if the minus sign holds and $(a-2)^2\le 1$ if the plus sign holds.
In the first case we get $a=0$ which finishes the proof.
So as only possibility $q=2$
$$\pmatrix{1&2&1\cr 2&4&2\cr 1&2&1}$$
remains. One can transform this matrix
by the matrix
$$\pmatrix{1&2&0\cr0&-1&0\cr0&0&-1}$$
to $E_6^-$. The other occurring matrices $E_2,E_3,E_4,E_5$ are fixed under this transformation.
This finishes the proof in the case A4. So case A is settled.
\medskip
It remains to treat the case B.
This case is very similar to the case A. Hence we can keep short.
Recall that the case B we consider the decomposition
$$T=r_1E_1+r_2E_2+r_3E_3+r_4E_4+r_5E_5+r_6E_6+r_7E_7$$
where
$$\eqalign{&r_1=t_{11}-t_{12}-t_{13}+m,\ r_2=t_{22}-t_{12}-t_{23}+m,
\ r_3=t_{33}-t_{13}-t_{23}+m,\cr
&r_4=t_{12}-m,\quad r_5=t_{23}-m,\quad r_6=t_{13}-m,\quad r_7=m.\cr}$$
The reduction conditions for $T$ imply that all $r_i$ are nonnegative.
At least one of the $r_4,r_5,r_6$ is zero. The remaining reduction inequalities are
$$\eqalign{&
r_1+r_6\le r_2+r_5,\quad r_2+r_4\le r_3+r_6,\quad r_4+r_7\le r_1+r_6,\cr
&r_5+r_7\le r_2+r_4,\quad r_6+r_7\le r_1+r_4,\cr& r_5+r_6+4r_7\le r_1+r_2.\cr}$$
Since the diagonal elements of $T$ are positive, we also have
$$r_1+r_4+r_6+r_7>0.$$
Again we differ between 4 cases where either all $r_1,r_2,r_3$ are positive or one of the is zero.
We claim that only the following 4 cases are possible,
\vskip1mm
\item{1)} $r_1>0$, $r_2>0$, $r_3>0$.
\item{2)} $r_1>0$, $r_2>0$, $r_3=0$ and $r_4=r_7=0$, $r_1=r_2=r_5=r_6$.
\item{3)} $r_1>0$, $r_2=0$, $r_3>0$ and $r_6=r_7=0$, $r_1=r_4=r_5$.
\item{4)}  $r_1=0$, $r_2>0$, $r_3>0$ and $r_6=r_4>0$, $r_5=r_7=0$.
\smallni
If $m$ is positive, then we are in the first case. Hence we can assume for the rest that $m=0$.
\smallskip
As  in the case A) we list the  corresponding matrices $T$. They
have the forms
$$ T= \pmatrix{2a&0&a \cr 0&2a&a\cr a&a&2a},\quad T= \pmatrix{2a&a&0 \cr a&2a&a\cr 0&a&c},\quad  T= \pmatrix{2a&a& a \cr a&b&0\cr  a&0&c} $$
 Really the case $B3)$   does not occur, since it contradicts $t_{13}>0$.
\smallni
{\it Proof of Proposition \ProOptz\ in the case B1.}
\smallni
As in the proof of A1 we can assume that $U_1=U_2=U_3=E$ and $U_i=D_i$ is diagonal for $i>3$
if $r_i\ne 0$
In the case $q=4$ the congruence $E_i[D_i]\equiv E_i$ mod 4 implies $E_i[D_i]=E_i$.
Hence we can assume $q=2$.
As we have shown during the proof of A1, we have $\lambda(E_4+E_5+E_6)>3$. Hence one of the
$r_4,r_5,r_6$ must be zero. The case $r_7=0$ is similar to the case A1 and can be omitted.
Hence we can assume $r_7>0$. There are three possibilities for $E_7[D_i]$ which behave similar.
We restrict to threat the case
$$E_7[D_i]=\pmatrix{1&1&-1\cr1&1&-1\cr -1&-1&1}.$$
Since $\lambda(E_1+E_5+E_7)=\lambda(E4+E_6^-+E_2+E_3)>3$ we must have $r_5=0$.
Similarly $\lambda(E_2+E_6+E_7)>3$ shows $E_6$. Now can apply the diagonal matrix with entries
$1,1,-1$. It transforms $E_7[D_7]$ to $E_7$ and keeps the other non-zero terms fixed.
\smallskip
In the cases B2)  and B4) the coefficient $r_7$ is zero. Hence we are in nearly the same
situation as in the cases  A2) and  A4).
This finishes the proof of Proposition~\ProOptz.\qed
\neupara{Localizations of rings of theta series}%
The algebra $R(g,q)$  of theta constants is generated by the theta constants
$$f_{a,q}=\sum_{n\in\gz^g}\exp\pii q Z[n+a/q].$$
We consider a decomposition $g=g_1+g_2$ and
$$Z=\pmatrix{\tau&z'\cr z&W},\quad \tau\in \calH_{g_1},\quad W\in\calH_{g_2}.$$
The Fourier expansion with respect to $W$ as variable can be written in the form
$$f_{a,q}=\sum_{T}f^T_{a,q}(\tau,z)\exp{\pii\over q}\tr (TW).$$
The coefficients $f^T_{a,q}$ can be considered as elements of $\calJ_T(\tau_0)$.
They can be different from 0 only if the rank of $T$ is $\le 1$.
\proclaim
{Proposition}
{Assume that $T$  is an integral semipositive $g_2\times g_2$-matrix of rank one
and with coprime entries, i.e.\ $T=a_2a'_2$, $a_2\in\gz^{g_2}$ coprime.
The $\calO_{\tau_0}$-module $\calJ_T(\tau_0)$ is generated by all
$$f^T_{a,q}(\tau,z)=\sum_{n_1\in\gz^{g_1}}\exp{\pii\over q}(\tau[qn_1+a_1]+
2a_2'z(qn_1+a_1)),\quad a=\pmatrix{a_1\cr a_2},$$
where $a_1$ runs through a system of representatives of $\gz^{g_1}/q\gz^{g_1}$.
}
ThNG%
\finishproclaim
{\it Proof.} The formula for the $f^T$ is obtained by a simple calculation.
By Nakayama's Lemma, it is sufficient to show that the vector space $J_T(\tau_0)$
is generated by the $f^T_{a,q}(\tau_0,z)$. They span a space of dimension
$q^{g_1}$. But this is also the dimension of this space (see [Ig4] for some explanations
about the dimensions of the spaces $J_T(\tau_0)$. Also the proof of Lemma \OZm\
can be extended to the computation of the dimensions.).\qed
\smallskip
The image of the cusp $\tau_0\in\calH_{g_1}$, $0\le g_1\le g$, in $\proj R(g,q)$
corresponds to the homogenous maximal ideal
$\gotm\subset R(g,q)$ consisting of all elements $f$ with the property
$$\lim_{t\to\infty} f\pmatrix{\tau_0&0\cr0&\imag tE}=0.$$
We consider its homogenous localization $R_{(\Gotm)}$.
It consists of quotients $f/g$, $g\not\in\gotm$, where $f,g$ are homogenous and of the same degree.
We are interested in cases where this
ring is normal.
\proclaim
{Lemma}
{The ring $R(g,q)_{(\Gotm)}$ is normal if and only if it is analytically irreducible
and if the ideal $\gotm$ generates the maximal ideal of $\hat S(\tau_0)$.}
AnIrr%
\finishproclaim
We recall that a local noetherian integral domain is analytically irreducible if its completion
is an integral domain. We denote by $\hat R(g,q)_{(\Gotm)}$ the completion of $R(g,q)_{(\Gotm)}$.
There is a natural homomorphism
$$\hat R(g,q)_{(\Gotm)}\lo \hat S(\tau_0).$$
It is surjective, since we assume that $\gotm$ generates the maximal ideal of $\hat S(\tau_0)$.
Since the left hand side is an integral domain by assumption, the map is an isomorphism.
Hence $\hat R(g,q)_{(\Gotm)}$ is a normal integral domain. This implies that
$R(g,q)_{(\Gotm)}$ is normal (by Zariski's main theorem).\qed
\smallskip
Igusa proved that in the case $q=4$ that the map $\overline{\calH_q/\Gamma_g[q,2q]}\to\proj R(g,q)$
is bijective. Therefore the local rings of the left hand side are analytically irreducible
in this case.
\proclaim
{Proposition}
{Assume that each $T\in\calT_{g_2}$ admits a $q$-optimal decomposition
$T=T_1+\cdots+T_k$ such that the multiplication map
$$J_{T_1}(\tau_0)\otimes_\cz\dots,\otimes_\cz J_{T_k}(\tau_0)\lo J_T(\tau_0)$$
is surjective. Then $\hat R(g,q)_{(\Gotm)}\lo \hat S(\tau_0)$ is surjective.}
AdSu
\finishproclaim
{\it Proof.} We have to show the following. Let $P$ be an element of the maximal ideal
of $\hat S(\tau_0)$. For each $k$ there exists an element $Q$ in the maximal ideal
of $R(g,q)_{(\Gotm)}$ such that $P-Q\in\hat{\gotm}_k$. It is sufficient to show
that for each $P\in\hat{\gotm}_k$ there exists $Q$ in the maximal ideal
of $R(g,q)_{(\Gotm)}$ such that $P-Q\in\hat{\gotm}_{k+1}$.
By definition of of $\hat{\gotm}_k$ we can write $P$ as a sum of products
$AB$ where $A$ is in $\gotm(\calO_{\tau_0})^\mu$ and where the coefficients of
$B$ are zero for $\lambda(T)<\nu$ and where $\mu+\nu=k$.
We can prove the statement separately for $A$ (with $\mu$ instead of $k$) and for
$B$ (with $\nu$ instead of $k$).
So it is  sufficient to assume that $P=A$ or $P=B$.
\smallni
{\it Case 1.} $P\in  \gotm(\calO_{\tau_0})^\mu$. In this case we can use
the result that the ring $R(g_1,q)$ gives a biholomorphic
embedding of $\calH_{g_1}/\Gamma_{g_1}[q,2q]$ into a projective space. Since the natural
projection $R(g,q)\to R(g_1,q)$ is surjective, this implies that the maximal
ideal of $\calO_{\tau_0}$ can be generated by (images of) linear combinations
of $f_{a,q}\in R(g,q)$ which vanish at $\tau_0$ divided by a suitable $f_{b,q}$ that does
not vanish at $\tau_0$.
\smallni
{\it Case 2.} The coefficients of
$P$ are zero for $\lambda(T)<k$. Then we choose an admissible decomposition
$T=T_1+\cdots+T_k$ and use the assumption in Proposition \AdSu.
This finishes the proof of this proposition.\qed
\smallskip
It remains to check whether the assumption of Proposition \AdSu\ is fulfilled.
We restrict now to $g=3$ and $q=4$. Then admissible decompositions exist.
We have to differ between three cases.
\smallni
{\it The case of a zero dimensional boundary component.} This case is trivial,
since in this case the spaces $J_T$ all are of dimension 1.
\smallni
{\it The case of a two dimensional boundary component.} In this case $T=m$ is a number.
The statement is that
$$J_1(\tau_0)^{\otimes m}\lo J_m(\tau_0)$$
is surjective. Since $J_1(\tau_0)$ is the space of sections of an ample line-bundle of the
form $\calL^4$, the statement follows from the well-known result that
$$H^0(\calL)^{\otimes m}\lo H^0(\calL^{\otimes m})$$
is surjective for $m\ge 3$.
\smallni
{\it The case of a one-dimensional boundary component.} The elements of $L_T(\tau_0)$
can be identified with the sections of a line-bundle on $E\times E$, where
$E=\cz/(\gz+\gz\tau_0)$.
In the case $T=E_1={1\,0\choose 0\,0}$ the space is spanned by the 4 theta series
$$\sum_{n\in\gz}\exp4\pii\{\tau(n+a_1/4)^2+2(n+a_1/4)z_1\}.$$
They can be considered as sections of a line bundle $\calL$ on the first component
$E$ of $E\times E$, i.e.~the line bundle on $E\times E$ is the inverse image
$\calL_1:=p^*\calL$ with respect to the first projection. Similarly in the case
$T=E_2={0\,0\choose 0\,1}$ we have to consider the line bundle $\calL_2:=q^*\calL$
where $q$ is the projection on the second $E$. Finally in the case
$T=E_3={1\,1\choose 1\,1}$ the line bundle $\calL_3=(p+q)^*\calL$ has to be considered.
\smallskip
We have to consider optimal decompositions of $2\times 2$-matrices. We can restrict to
reduced matrices, then the optimal decompositions are of the form
$$T=aE_1+bE_2+cE_3.$$
What we have to show is that the multiplication map
$$H^0(\calL_1)^{\otimes a}\otimes_\cz H^0(\calL_2)^{\otimes b}\otimes_\cz H^0(\calL_3)^{\otimes c}
\lo H^0(\calL_1^{\otimes a}\calL_2^{\otimes b}\calL_3^{\otimes c})$$
is surjective. This is the problem for the cartesian square of an elliptic curve.
All what we must know is that $\calL$  is a line-bundle (= divisor class) on $E$ of degree
4. But this follows from $\dim H^0(\calL)=4$. Any  divisor of degree 4 is equivalent to
a translate of $4[0]$. Since we are free to change the origin we can assume that
$\calL$ is the line bundle associated to the divisor $4[0]$.
So we can reformulate the problem as follows.
\proclaim
{Proposition}
{We denote by $L(a,b,c)$ the space of all meromorphic functions on $E\times  E$
which are regular or have poles of order $\le a$ on $\{0\}\times E$, of order $\le b$
on $E\times \{0\}$ and of order $\le c$ on the diagonal. The multiplication map
$$L(4,0,0)^{\otimes a}\otimes_\cz L(0,4,0)^{\otimes b}\otimes_\cz L(0,0,4)^{\otimes c}
\lo L(4a,4b,4c)$$
is surjective.}
Labc%
\finishproclaim
We shall prove this proposition in the next section.
It will include our main result.
\proclaim
{Main-Theorem}
{In the case $g=3$ the theta functions $f_a$, $a\in(\gz/4\gz)^6/\pm$, define a biholomorphic
embedding of the Satake compactification $\overline{\calH_3/\Gamma_3[4,8]}$ into the projective
space.}
MT%
\finishproclaim
As we mentioned already one can replace the $f_a$ by the standard 36 theta constants of first
kind.
\neupara{Cartesian square of an elliptic curve}%
In this section we give the proof of Proposition \Labc\ (and hence of Main-Theorem~\MT).
We consider the elliptic curve $E=\cz(\gz+\gz\tau)$, $\im\tau>0$.
We will construct the spaces $L(a,b,c)$ (see Proposition \Labc) by means
of the Weierstrass $\wp$-function.
We will use the basic fact that every elliptic function (meromorphic function on $E$)
that is holomorphic outside the origin can be written as unique linear combination of the
(higher) derivatives of $\wp$ including $\wp$ and the constant function 1.
We consider the matrix
\smallni
\halign{$\qquad#$\quad\hfill&$#$\quad\hfill&$#$\quad\hfill&$#$\quad\hfill&$#$\quad\hfill\cr
1,&\wp(z),&\wp'(z),&\wp''(z),\cr
1,&\wp(w),&\wp'(w),&\wp''(w),\cr
1,&\wp(z-w),&\wp'(z-w),&\wp''(z-w).\cr
}
\smallni
If we take from the first line $a$ elements, from the second line $b$
elements and from the third line $c$ elements and multiply them, we get three-products
in $L(4a,4b,4c)$ .
We denote the subspace of $L(4a,4b,4c)$ generated by them by $M(4a,4b,4c)$.
So the statement of Proposition \Labc\ is $L(4a,4b,4c)=M(4a,4b,4c)$.
We notice that
$$\wp^{(k)}(z)\in M(4a,0,0),\quad\hbox{if}\ k+2\le 4a.$$
\proclaim
{Lemma}
{The function
$$\varphi(z,w)={\wp'(z)+\wp'(w)\over \wp(z)-\wp(w)}$$
has poles of first order at
the $3$ special divisors
($E\times\{0\}$, $\{0\}\times E$ and diagonal)
and has no other pole.
Hence it is contained in $L(1,1,1)$.}
PhiM%
\finishproclaim
{\it Proof.\/}
This follows form the addition formula for the $\wp$-function,
$$\varphi(z,w)^2=\wp(z-w)+\wp(z)+\wp(w).\eqno\square$$
\par We consider the matrix group $G$ generated by the matrices
$$\pmatrix{0&1\cr 1&0},\quad \pmatrix{1&-1\cr0&-1}.$$
It has order 6.
It contains the negative unit matrix in its center.
It acts on the variables (z,w) trough
$$\pmatrix{z\cr w}\loma g \pmatrix{z\cr w}$$
from the left and hence on functions in $z,w$ from the right.
\proclaim
{Lemma}
{The function $\varphi(z,w)$ has the property
$$\varphi(z,w)=\varphi(-z,-w).$$
Moreover, it is invariant under $G$ up to the character
$\varepsilon(g)=\det(g)$. The formula
$$\lim_{z\to 0}z\varphi(z,w)=-2$$
holds.
}
PhiCh%
\finishproclaim
We want to exhibit all functions from
the space $M(4,4,4)$ that have the same transformation formula as $\varphi$.
The 64 generating functions all are symmetric or skew symmetric under
$(z,w)\mapsto (-z,-w)$. We are only interested in the skew symmetric ones.
The group $G$ acts on them. A system of representatives is given by
the functions
$$\eqalign{&
\wp'(z),\quad \wp(z)\wp'(w),\quad \wp''(z)\wp'(w),\quad
\wp(z)\wp(w)\wp'(z-w),\cr&
\wp(z)\wp''(w)\wp'(z-w),\quad \wp''(z)\wp''(w)\wp'(z-w),\quad \wp'(z)\wp'(w)\wp'(z-w).\cr}$$
We symmetrize them with respect to the Character $\varepsilon$.
$$\eqalign{
&f_1=\ \wp'(z)-\wp'(w)-\wp'(z-w),\cr
\noalign{\vskip1mm}
&f_2=\ \wp(z)\wp'(w)-\wp'(z)\wp(w)+\wp'(w)\wp(z-w)+\wp(w)\wp'(z-w)-\cr
&\qquad\quad\wp'(z)\wp(z-w)+\wp(z)\wp'(z-w),\cr
\noalign{\vskip1mm}
&f_3=\ \wp''(z)\wp'(w)-\wp'(z)\wp''(w)+\wp'(w)\wp''(z-w)+\wp''(w)\wp'(z-w)-\cr
&\qquad\quad\wp'(z)\wp''(z-w)+\wp''(z)\wp'(z-w),\cr
\noalign{\vskip1mm}
&f_4=\ \wp(z)\wp(w)\wp'(z-w)-\wp'(z)\wp(w)\wp(z-w)+\wp(z)\wp'(w)\wp(z-w),\cr
\noalign{\vskip1mm}
&f_5=\ \wp(z)\wp''(w)\wp'(z-w)+\wp''(z)\wp(w)\wp'(z-w)-\wp'(z)\wp''(w)\wp(z-w)-\cr
&\qquad\quad\wp'(z)\wp(w)\wp''(z-w)+\wp''(z)\wp'(w)\wp(z-w)+\wp(z)\wp'(w)\wp''(z-w),\cr
\noalign{\vskip1mm}
&f_6=\ \wp''(z)\wp''(w)\wp'(z-w)-\wp'(z)\wp''(w)\wp''(z-w)+\wp''(z)\wp'(w)\wp''(z-w),\cr
&f_7=\ \wp'(z)\wp'(w)\wp'(z-w).\cr
}$$
We have to investigate the pole behavior of these functions. The only poles are along the three divisors.
The symmetry properties show that the behavior at each of the three divisors is the same.
Hence it is sufficient to concentrate on the divisor $z=0$. We compute some Laurent coefficients
for fixed $w=a\ne 0$. What we need is
$$\wp(z)={1\over z^2}+O(1),\quad \wp'(z)=-{2\over z^3}+O(1),\quad \wp''(z)={6\over z^4}+O(1).$$
Here $O(1)$ stands for a bounded function in a small neighborhood of the origin.
We also need
$$\eqalign{
&\wp(z-a)=\wp(a)-\wp'(a)z+{\wp''(a)\over 2}z^2-{\wp^{(3)}(a)\over 6}z^3+\cdots,\cr
&\wp'(z-a)=-\wp'(a)+\wp''(a)z-{\wp^{(3)}(a)\over 2}z^2+{\wp^{(4)}\over 6}z^3+\cdots,\cr
&\wp''(z-a)=\wp''(a)-\wp^{(3)}(a)z+{\wp^{(4)}\over 2}z^2-{\wp^{(5)}\over 6}z^3+\cdots,\cr}$$
By means of these formualae we are able to compute the Laurent coefficients. The differential equation
of the $\wp$-function allows the express the higher derivatives of $\wp$ explicitly in terms
of $\wp$ and $\wp'$. Now a somewhat tedious but straightforward computation gives the
following result.
\proclaim
{Proposition}
{The functions
$$\eqalign{F_1&=f_3-30f_4-(5/2)g_2f_1\cr
F_2&=2f_5-9f_7-5g_2f_2+15f_1\cr
}$$
have poles of order
$1$ along the three divisors $D_i$. They are contained
in $L(1,1,1)\cap M(4,4,4)$.
We have
$$\eqalign{
\lim_{z\to 0}zF_1(z,w)&=-36g_2\wp(w)-54g_3,\cr
\lim_{z\to 0}zF_2(z,w)&=108g_3\wp(w)+6g_2^2
.}$$
The functions $1$, $F_1$ and $F_2$ are linearly independent. They span the space $L(1,1,1)$.
}
FeFz%
\finishproclaim
In particular, $\varphi$ must be a linear combination of $F_1$ and $F_2$. Here is it.
\proclaim
{Proposition}
{We have
$$3(g_2^3-27g_3^2)\varphi=-g_2F_2+3g_3F_1.$$}
PhiF%
\finishproclaim
Notice that the discrimant $g_2^3-27g_3^2$ is different from 0.
\smallni
\proclaim
{Proposition}
{The spaces $L(4a,4b,4c)$ and $M(4a,4b,4c)$ agree.}
WpAll%
\finishproclaim
{\it Proof.} Let $f(z,w)\in L(4a,4b,4c)$.
In  a first step we assume that the order of $f$ along one of the three components
is zero. Without loss of generality we can assume that the order at the diagonal is
zero. Then, for fixed $w\ne 0$ the function $z\mapsto f(z,w)$ has only a pole at $z=0$.
Hence it can be written as linear combination  in the derivatives of the $\wp$-function
(including the constant function),
$$f(z,w)=a_0+\sum_{\nu\ge0}a_\nu\wp^{(\nu)}(z).$$
The coefficients $a_\nu$ are elliptic functions in $w$ with poles only at $w=0$. Hence they
can be expressed by derivatives of $\wp(w)$ (including the function constant 1 and $\wp(w)$).
\smallskip
Now we can assume that the order of $f$ along $z=0$ is $m>0$. We first treat the case where
$m>1$. Again we fix $w\ne 0$. Then $f(z,w)$ has a pole of order $>1$ at $z=0$ (and may be a pole
at $z=w$). We subtract from $f(z,w)$ a constant multiple of $\wp^{(m-2)}(z)$ such that the difference
$f(z,w)-a \wp^{(k-2)}(z)$ has smaller pole order at $z=0$. Again the coefficient $a=a(w)$
is an elliptic function with poles only at $w=0$. It can be expressed by derivatives of
$\wp(w)$.
\smallskip
In the remaining case $m=1$ we consider for fixed $w$ the difference
$$f(z,w)-a\varphi(z,w).$$
The pole at $z=0$ can be cancelled.
The coefficient $a$ is an elliptic function in $w$ with now poles outside $w=0$. Hence it can be
expressed by derivatives of $\wp(w)$.
This finishes the proof of Proposition \Labc\ and hence of our main result.\qed
\smallskip
We mention that In Proposition \WpAll\ the factor 4 is essential. For example
$L(1,1,1)$ is not contained in $M(3,3,3)$.
\vskip1cm\noindent
\paragratit References\rm
\medni
\item{[FK]} Freitag, E., Kiehl, R.:
{\it Algebraische Eigenschaften der lokalen Ringe in den Spitzen der Hilbertschen Modulgruppe,}
Inv.Math. {\bf 24}, 121-148 (1974)
\medni
\item{[Ig1]} Igusa, J.I.: {\it On the Graded Ring of Theta-Constants,}
Am.\ J.\ of Math.
Vol.\ {\bf 86}, No.\ 1, pp.\ 219-246 (1964)
\medskip
\item{[Ig2]} Igusa, J.I.: {\it On the Graded Ring of Theta-Constants II,}
Am.\ J.\ of Math.
Vol.\ {\bf 88}, No.\ 1, pp.\ 221-236 (1966)
\medskip
\item{[Ig3]} Igusa, J.I.: {\it On the Variety Associated with the Ring of Thetanullwerte,}
Am.\ J.\ of Math.
Vol.\ {\bf 103}, No.\ 2, pp. 377--398, (1981)
\medskip
\item{[Ig4]} Igusa, J.I.: {\it A desingularization problem in the theory of Siegel
modular varieties,} Math.\ Ann. Vol.\ {\bf 168}, 228--260  (1967)
\medskip
\item{[Ig5]}  Igusa, J.I.: {\it On the nullwerte of Jacobians of odd theta
functions,} Symposia Mathematica, Vol. XXIV (Sympos., INDAM, Rome,
1979), pp.\ 83–-95, Academic Press, London-New York (1981)
\medskip
\item{[Kn]} Kn\"oller, F.W.: {\it Multiplizit\"aten
``unendlich ferner'' Spitzen,} Mh.\ Math. {\bf 88}, 7--26 (1979)
\medskip
\item{[Mu]} Mumford, D., with Nori, M. and Norman, P.: {\it Tata Lectures on Theta III,}
Modern Birkhäuser Classics (2010)
\medskip
\item{[Ru]} Runge, B. {\it On Siegel modular forms. I}
J. Reine Angew. Math. 436, 57-85 (1993).
\medskip
 \item{[SM]} Salvati Manni, R. {\it On the projective varieties associated with some subrings of the ring of Thetanullwerte}
Nagoya Math. J. 133, 71-83 (1994).
\bye